\documentclass[12pt,reqno]{amsart}
\usepackage{amsmath , amssymb }
\usepackage{graphics}

\usepackage[mathscr]{eucal}

\newdimen\unit\newdimen\psep\newcount\nd\newcount\ndx\newbox\dotb\newbox\ptbox
\newdimen\dx\newdimen\dy\newdimen\dxx\newdimen\dyy\newdimen\hgt
\newdimen\xoff\newdimen\yoff
\newcommand\clap[1]{\hbox to 0pt{\hss{#1}\hss}}
\newcommand\vdisk[1]{{\font\dotf=cmr10 scaled #1\dotf.}}
\newcommand\varline[2]{\setbox\dotb\hbox{\vdisk{#1}}\xoff=-.5\wd\dotb
\wd\dotb=0pt\yoff=-.5\ht\dotb\psep=#2\ht\dotb}
\newcommand\varpt[1]{\setbox\ptbox\clap{\vdisk{#1}}\setbox\ptbox
\hbox{\raise-.5\ht\ptbox\box\ptbox}}
\newcommand\cpt{\copy\ptbox}
\newcommand\point[3]{\rlap{\kern#1\unit\raise#2\unit\hbox{#3}}}
\newcommand\setnd[4]{\dx=#3\unit\advance\dx-#1\unit\divide\dx by\psep
\dy=#4\unit\advance\dy-#2\unit\divide\dy by\psep \multiply\dx
by\dx\multiply\dy by\dy\advance\dx\dy\nd=1\advance\dx-1sp
\loop\ifnum\dx>0\advance\dx-\nd sp\advance\nd1\advance\dx-\nd
sp\repeat}

\newcommand\dline[5]{{\nd=#5\hgt=#2\unit\dx=#3\unit\advance\dx-#1\unit
\divide\dx by\nd\dy=#4\unit\advance\dy-#2\unit\divide\dy by\nd
\advance\hgt\yoff\rlap{\kern#1\unit\kern\xoff\loop\ifnum\nd>1\advance\nd-1
\advance\hgt\dy\kern\dx\raise\hgt\copy\dotb\repeat}}}

\newcommand\qellip[4]{{\setnd{0}{0}{#3}{#4}\dx=\unit\dy=0pt\raise\yoff\rlap{%
\kern#1\unit\kern\xoff\raise#2\unit\hbox{\loop\ifnum\dx>0\rlap{\kern#3\dx
\raise#4\dy\copy\dotb}\hgt=\dx\divide\hgt
by\nd\advance\dy\hgt\hgt=\dy \divide\hgt
by\nd\advance\dx-\hgt\repeat\rlap{\raise#4\dy\copy\dotb}}}}}
\newcommand\bez[6]{{\setnd{#1}{#2}{#3}{#4}\ndx=\nd\setnd{#3}{#4}{#5}{#6}
\ifnum\ndx>\nd\nd=\ndx\fi\dx=#3\unit\advance\dx-#1\unit\dy=#4\unit
\advance\dy-#2\unit\dxx=#5\unit\advance\dxx-#1\unit\dyy=#6\unit\advance
\dyy-#2\unit\advance\dxx-2\dx\advance\dyy-2\dy\divide\dxx
by\nd\divide\dyy
by\nd\advance\dx.25\dxx\advance\dy.25\dyy\divide\dx
by\nd\divide\dy by\nd \multiply\nd
by2\dx=100\dx\dy=100\dy\dxx=100\dxx\dyy=100\dyy\divide\dxx by\nd
\divide\dyy
by\nd\hgt=#2\unit\raise\yoff\rlap{\kern#1\unit\kern\xoff
\raise\hgt\copy\dotb\loop\ifnum\nd>0\advance\nd-1\advance\hgt0.01\dy
\kern0.01\dx\raise\hgt\copy\dotb\advance\dx\dxx\advance\dy\dyy\repeat}}}
\newcommand\ptu[3]{\point{#1}{#2}{\cpt\raise1ex\clap{$\scriptstyle{#3}$}}}
\newcommand\ptd[3]{\point{#1}{#2}{\cpt\raise-1.8ex\clap{$\scriptstyle{#3}$}}}
\newcommand\ptr[3]{\point{#1}{#2}{\cpt\raise-.4ex\rlap{$\ \scriptstyle{#3}$}}}
\newcommand\ptl[3]{\point{#1}{#2}{\cpt\raise-.4ex\llap{$\scriptstyle{#3}\ $}}}
\newcommand\ptlu[3]{\point{#1}{#2}{\raise.8ex\clap{$\scriptstyle{#3}$}}}
\newcommand\ptld[3]{\point{#1}{#2}{\raise-1.6ex\clap{$\scriptstyle{#3}$}}}
\newcommand\ptlr[3]{\point{#1}{#2}{\raise-.4ex\rlap{$\,\scriptstyle{#3}$}}}
\newcommand\ptll[3]{\point{#1}{#2}{\raise-.4ex\llap{$\scriptstyle{#3}\,$}}}

\newcommand\thnline{\varline{400}{.6}}

\varpt{2500}\thnline\unit=2em

\newtheorem{thm}{Theorem}
                                                                                   
\newtheorem{conj}{Conjecture}
\newtheorem{qu}{Question}
\newtheorem{prob}{Problem}
\newtheorem{lemma}[thm]{Lemma}

\newtheorem{counter}{Counter-example}
\newtheorem{obs}[thm]{Observation}
\theoremstyle{definition}\newtheorem{rmk}{Remark}
\theoremstyle{definition}
\newtheorem*{conjFC}{The Union-Closed Sets Conjecture}
\newtheorem*{thmP}{Poonen's Theorem}
\newtheorem*{lemmaA}{Lemma A}
\newtheorem*{lemmaB}{Lemma B}
\newtheorem*{thmC}{Theorem C}
\newtheorem*{thmD}{Theorem D}
\newtheorem*{lemmaE}{Lemma E}
\newcommand{\ds}{\displaystyle}
\newcommand{\ul}{\underline}

\def\A{\mathcal{A}}
\def\B{\mathcal{B}}

\def\P{\mathcal{P}}

\def\S{\mathcal{S}}

\def\le{\leqslant}
\def\ge{\geqslant}

\def\eps{\varepsilon}

\begin{document}
\title[FC-families and Frankl's Conjecture]{FC-families, and improved bounds for Frankl's Conjecture}

\author{Robert Morris}
\address{Department of Mathematical Sciences\\ The University of Memphis\\ Memphis, TN 38152} \email{rdmorrs1@memphis.edu}\thanks{The author was supported during this research by a Van Vleet Memorial Doctoral Fellowship.}

\begin{abstract}
A family of sets $\A$ is said to be union-closed if $\{A \cup B : A,B \in \A\} \subset \A$. Frankl's conjecture states that given any finite union-closed family of sets, not all empty, there exists an element contained in at least half of the sets. Here we prove that the conjecture holds for families containing three 3-subsets of a 5-set, four 3-subsets of a 6-set, or eight 4-subsets of a 6-set, extending work of Poonen and Vaughan. As an application we prove the conjecture in the case that the largest set has at most nine elements, extending a result of Gao and Yu. We also pose several open questions.
\end{abstract}

\maketitle

\section{Introduction}\label{intro}

A family of sets $\A$ is called \textit{union-closed} if $A \cup B \in \A$ for every pair of sets $A, B \in \A$. The Union-Closed Sets Conjecture, which is also known as Frankl's Conjecture, and is generally attributed to P.~Frankl in 1979 (see \cite{Duff}, \cite{Fran}), is as follows.

\begin{conjFC}
If $\A$ is a finite union-closed family of sets, not all empty, then there exists an element belonging to at least half of the sets of $\A$.
\end{conjFC}

\noindent By Corollary 1 of \cite{Poon} we may assume that the sets of $\A$ are finite, so, writing $\P(n)$ for the power set of $\{1, \ldots, n\}$, let $\A \subset \P(n)$ and $|\A| = m$.

Very little progress has been made on this simple-sounding problem. David Reimer~\cite{Reim} recently proved that the average set-size in a union-closed family is at least $\ds\frac{1}{2}\log_2(m)$, but this gives only weak bounds on the number of appearances of the most common element. We shall continue a route into the problem initiated by Poonen~\cite{Poon}, and continued by Vaughan~\cite{Vaug}, \cite{Vnote}, \cite{Vthrees}, and also by Gao and Yu~\cite{GY}. The approach is motivated by the following simple observations. First, let $\A$ be a union-closed family of sets, and suppose that $\A$ contains a one-element set $\{x\}$. Then clearly $x$ is in at least half the sets of $\A$. Similarly, if $\A$ contains a two-element set $\{x,y\}$, then at least one of $x$ and $y$ is contained in at least half the sets of $\A$ (consider the four sub-families containing both, exactly one, and neither of the pair $x$ and $y$).

Following Vaughan~\cite{Vaug}, call $\B$ an \textit{$FC(k)$-family} if it is a union-closed family of sets containing $\emptyset$ and whose largest set has $k$ elements, and it has the property that given any union-closed family $\A \supseteq \B$, one of the $k$ elements of (the largest set in) $\B$ is in at least half the sets of $\A$. Call a family $FC$ if it is $FC(k)$ for some $k$, and call an $FC(k)$-family \textit{proper} if it contains no strictly smaller $FC$-family. By the observations above, $\{\emptyset, \{x\}\}$ is an $FC(1)$-family, and $\{ \emptyset, \{x,y\} \}$ is an $FC(2)$-family. However, a single 3-element set is not an $FC(3)$-family (an example showing this was given in \cite{SR}), so there are no proper $FC(3)$-families.

Poonen studied these families, and gave necessary and sufficient conditions for a family to be $FC$ (see Poonen's Theorem, below). He used his result to show that the conjecture holds for any family which contains three of the 3-subsets of a 4-set, and deduced that the conjecture holds if $n \le 7$ or $m \le 28$. Gao and Yu~\cite{GY} later improved these bounds to $n \le 8$ and $m \le 32$. Vaughan~\cite{Vaug}, \cite{Vnote}, \cite{Vthrees} studied $FC(k)$-families for small values of $k$, and showed that the conjecture holds for any family which contains all five of the 4-subsets of a 5-set, or ten of the 4-subsets of a 6-set, or three 3-subsets of a 7-set with a common element.
Given a set system $\S \subset \P(n)$, we define the \textit{family generated by} $\S$ to be the smallest union-closed family containing $S \cup \emptyset$, i.e. $\{A \in \P(n) : A = A_1 \cup \ldots \cup A_r$, \; where $A_i \in \S \cup \emptyset$ \: for \: $1 \le i \le r \}$. Since determining exactly which families are $FC$ seems to be complicated for $k \ge 5$, we shall concentrate on a slightly simpler question: how many $k$-sets in $[n] = \{1, \ldots ,n\}$ necessarily generate an $FC$-family? To this end, write $FC(k,n)$ for the minimal $m$ such that any $m$ of the $k$-sets in [$n$] generate an $FC$-family.

As noted above, Poonen~\cite{Poon} showed that $FC(3,4) = 3$, and Vaughan~\cite{Vaug} showed that $FC(4,5) \le 5$ and $FC(4,6) \le 10$. In this paper we shall improve on the results of Vaughan, by determining exactly which families are $FC(5)$, and finding new examples of $FC(6)$-families. In particular, we shall show that $FC(3,5) = 3$, $FC(4,5) = 5$, $FC(3,6) = 4$ and $7 \le FC(4,6) \le 8$. In other words, we shall prove that the conjecture holds for any union-closed family which contains any three of the 3-subsets of some 5-set, four of the 3-subsets of some 6-set, or eight of the 4-subsets of some 6-set. As an application of our results, we shall prove the conjecture in the cases $n \le 9$ and $m \le 36$, improving the bounds of Gao and Yu. We shall also prove a conjecture of Vaughan~\cite{Vaug} on union-closed families with many 2-sets, and discuss the asymptotics of $FC(k,n)$ as $n \to \infty$.

In this paper we shall prove the following theorems. The first determines exactly which families are $FC(5)$.

\begin{thm}\label{fvs}
A sub-family of $\P(5)$ is a proper \textup{$FC(5)$}-family if, and only if it is generated by one of the following set systems (under some permutation of $\{1,2,3,4,5\}$).
\begin{enumerate}
\item[$(a)$] Any three of the $3$-subsets,
\item[$(b)$] $\{1,2,3\}$, $\{1,2,4\}$ and $\{1,3,4,5\}$,
\item[$(c)$] $\{1,2,3\}$, $\{1,4,5\}$ and $\{2,3,4,5\}$,
\item[$(d)$] $\{1,2,3\}$, $\{1,4,5\}$, $\{1,2,3,4\}$, $\{1,2,3,5\}$, $\{1,2,4,5\}$ and $\{1,3,4,5\}$,
\item[$(e)$] $\{1,2,3\}$, $\{1,2,4,5\}$, $\{1,3,4,5\}$ and $\{2,3,4,5\}$,
\item[$(f)$] All five $4$-sets.
\end{enumerate}
In particular $FC(3,5) = 3$ and $FC(4,5) = 5$.
\end{thm}

The next gives some new examples of $FC(6)$-families.

\begin{thm}\label{sxs}
$FC(3,6) = 4$ and \hspace{0.01cm} $7 \le FC(4,6) \le 8$.
\end{thm}

As an application of our results, we prove the following theorem.

\begin{thm}\label{nine}
The Union-Closed Sets Conjecture holds in the case $n = 9$.
\end{thm}

We shall also prove the following theorem on union-closed families, conjectured by Vaughan~\cite{Vaug}, and proved by her in the cases $n = 5, 6$ and $7$, which will be useful in proving the above results.

\begin{thm}\label{vcj}
If $\A$ is a union-closed family in $\P(n)$ and contains at least $\ds{{n-1} \choose 2}+1$ of the $2$-sets, then
$$|\{A\in \A : |A| = n - k \}| \ge |\{A\in \A : |A| = k \}|$$ for every $0 \le k \le \left\lfloor \ds\frac{n}{2} \right\rfloor$.
\end{thm}

The rest of the paper is organised as follows. In Section~\ref{methodsec} we shall describe the tools we shall use to prove our main theorems; in Section~\ref{Vpfsec} we shall prove Theorem~\ref{vcj}; in Section~\ref{smallsec} we shall prove Theorems~\ref{fvs} and \ref{sxs}; in Section~\ref{9sec} we shall prove Theorem~\ref{nine}; in Section~\ref{bigsec} we shall discuss the asymptotics of $FC(k,n)$; and in Section~\ref{QCsec} we shall discuss some of our questions and conjectures relating to this work.

\section{Our method of approach}\label{methodsec}

Our main tool in finding $FC$-families will be the following theorem of Poonen~\cite{Poon}.

\begin{thmP}\label{poonen}
If $\B$ is a union-closed family containing $\emptyset$ whose largest set has $n$ elements, say $\{1, \ldots ,n\}$, then the following are equivalent:
\begin{enumerate}
\item[$(a)$] $\B$ is an $FC(n)$-family
\item[$(b)$] There exist non-negative real numbers $c_1, \ldots ,c_n$ summing to $1$ such that for every union-closed family $\A \subset \P(n)$ such that $\A \uplus \B \subseteq \A$,
$$\ds\sum_{i = 1}^n c_i|A_i|\ge \ds\frac{|A|}{2}$$
where $\A \uplus \B = \{A \cup B : A \in \A, B \in \B \}$ and $A_i = \{A \in \A : i \in A \}$.
\end{enumerate}
\end{thmP}

\noindent It will be convenient to allow the $c_i$ to be integers, and to define
$$K(\A) = \sum_{A \in \A} \left( \sum_{i \in A} c_i - \sum_{i \notin A} c_i \right) = 2\sum_{i=1}^n \sum_{A \in \A} \left( c_iI[i \in A] \right) - |\A| \sum_{i = 1}^n c_i,$$ where $I[S]$ as usual denotes the indicator function of the event $S$. Write $N_i$ for $\ds\sum_{\{A \in \A : |A| = i\}} \left( \ds\sum_{j \in A} c_j - \ds\sum_{j \notin A} c_j \right)$, the contribution of the $i$-sets to $K(\A)$, and notice that condition 2 of the theorem holds for a given family $A$ if and only if $K(\A) = \ds\sum_{i=0}^n N_i \ge 0$. Whenever possible, we shall choose the $c_i$ to be integers, but the reader should be aware that when proving that no such $c_i$ exist (for certain families), we shall revert to real numbers summing to one. It will always be clear which situation we are in, and we trust this will not cause any confusion. We shall write $\ul{c} = 1$ if $c_i = 1$ for all $i$, and $n_i$ for the number of $i$-sets in $\A$.

Our proof will also use the following lemmas of Vaughan~\cite{Vaug}, and the trivial observation below.

\begin{lemmaA}\label{djk}
Suppose that $\A$ and $\B$ are union-closed families in $\P(n)$ such that $\A \uplus \B \subset \A$, and that there are exactly $r$ $k$-sets which are in $\B$ but not $\A$. Suppose also that $r = \ds{{a_k} \choose k} + {{a_{k-1}} \choose {k-1}} + \ldots + {{a_t} \choose t}$ with $a_k > \ldots > a_t \ge t \ge 1$.

Then for each $1 \le j \le k$, the number of $j$-sets in $\A$ is bounded above by $d_{j}(n,k,r)$, where $$d_{j}(n,k,r) = \ds{n \choose j} - {{a_k} \choose j} - {{a_{k-1}} \choose {j-1}} - \ldots - {{a_t} \choose {t-k+j}}$$
\end{lemmaA}

\begin{lemmaB}\label{ti}
Let $t(i) = |\{A \in \A : i \in A, |A| = 2 \}|$. If $t(i) \ge k + 1$ for all $i \in [n]$, then $\A$ contains all of the $(n - k)$-sets in $\P(n)$.
\end{lemmaB}

\begin{rmk}
Note that in general if $t_j(i) = |\{A \in \A : i \in A, |A| = j \}|$ and $t_j(i) \ge \ds{{n-1} \choose {j-1}} - {{n-k-1} \choose {j-1}} + 1$ for some $j$ and all $i \in [n]$, then $\A$ contains all of the $(n - k)$-sets in $\P(n)$, since for any $(n-k)$-set $K$ and any $i \in K$ there exists a $j$-set $L$ such that $i \in L \subset K$.
\end{rmk}

\begin{obs}\label{triv}
If $\A$ and $\B$ are union-closed families with $\A \uplus \B \subseteq \A$, and $\emptyset \in \A$, then $\B \subseteq \A$.
\end{obs}

\section{Proof of Theorem~\ref{vcj}}\label{Vpfsec}

\begin{proof}
The cases $k = 0$ and $k = \ds\frac{n}{2}$ are trivial, so fix $1 \le k \le \left\lfloor \ds\frac{n-1}{2} \right\rfloor$ and suppose that the theorem fails for $k$ (i.e. $\A$ contains at least $\ds{{n-1} \choose 2}+1$ of the $2$-sets and more $k$-sets than ($n-k$)-sets). Note that we may assume that $\A$ is generated by its 2-sets and its $k$-sets, since the addition of other sets can only increase the number of ($n-k$)-sets. Consider the graph $G$ with [$n$] as the vertices and $ij \in E(G) \Leftrightarrow \{i,j\}\in \A$ (so subgraphs of $G$ with no isolated vertices correspond to sets in $\A$), and note that if each vertex of $G$ has degree at least $k + 1$ then we are done by Lemma B.\\

\noindent\ul{Claim}: Suppose either more than one element has degree at most $k$, or exactly one element has this property and in the induced subgraph on the other $n-1$ vertices some vertex has degree at most $k$. Then $n$ is odd, $k = \ds\frac{n-1}{2}$ and $G$ contains an induced copy of $K_{n-2}$.

\begin{proof}[Proof of claim]
We must be missing at least $n-1-k$ edges from the first vertex, and at least $n-2-k$ different edges from the second vertex. However, at most $n-2$ edges are missing from $G$, thus $n - 2 \ge (n-1-k) + (n-2-k)$, and so $k \ge \ds\frac{n-1}{2}$. It follows that $k = \ds\frac{n-1}{2}$, hence $n$ is odd, and equality holds everywhere, so all missing edges are incident with one of these two vertices.
\end{proof}

First assume that the assumptions of the claim fail to hold. Thus we may assume that exactly one vertex has degree at most $k$. Let this element be $u$, suppose $u$ lies in exactly $r$ of the $2$-sets in $\A$ (so $1 \le r \le k \le \left\lfloor \ds\frac{n-1}{2} \right\rfloor$). Note also that, by Lemma B, $\A$ contains all of the ($n-k$)-sets and all of the ($n-k-1$)-sets in $\P([n]-u)$.

Consider which of the ($n-k$)-sets in $\P(n)$ may be missing from $\A$. Such a set must contain $u$ by the comment above, and cannot contain any element of $\Gamma(u) = \{v : \{u,v\}\in \A \}$, as then it would be the union of some ($n-k-1$)-set in $\P([n]-u)$ with some pair $\{u,v\}\in \A$. So there are at most $\ds{{n-r-1} \choose {n-k-1}}$ such sets.

Now observe that any $k$-set which contains $u$ but no neighbour of $u$ is not generated by the 2-sets of $\A$. There are $\ds{{n-r-1} \choose {k-1}}$ such sets. Since $\ds{{n-r-1} \choose {n-k-1}} \le {{n-r-1} \choose {k-1}}$ when $r \ge 1$, we must have some of these sets in $\A$ as generating sets. However, if such a $k$-set $A\in \A$, then every ($n-k$)-set containing $A$ and $n-2k$ different elements of $\P([n]-u-\Gamma(u))$ is also in $\A$, by taking the union of $A$ with some ($n-k-1$)-set in $\P([n]-u-\Gamma(u))$ containing $A \setminus u$. Moreover each of these $(n-k)$-sets was counted as missing from $\A$ above.

In this way, each of these $k$-sets generates $\ds{{n-k-r} \choose {n-2k}}$ new ($n-k$)-sets, and each ($n-k$)-set is generated by $\ds{{n-k-1} \choose {k-1}}$ of the $k$-sets. But now simply observe that $\ds{{n-k-r} \choose {n-2k}} = {{n-k-r} \choose {k-r}} \le {{n-k-1} \choose {k-1}}$ when $r \ge 1$, and it is clear that the number of missing ($n-k$)-sets must be no greater than the number of missing $k$-sets.

We are left to deal with the (easier) case that $n$ is odd, $k = \ds\frac{n-1}{2}$ and $G$ contains an induced copy of $K_{n-2}$. Let the remaining two vertices be $u$ and $v$ and let $d(u) \le d(v)$. We may assume that the assumptions of the claim hold, and that $e(G) = \ds{{n-1} \choose 2} + 1$.

Suppose first that $uv \in E(G)$. Then $u$ and $v$ must be missing $\ds\frac{n-1}{2}$ and $\ds\frac{n-3}{2}$ edges respectively, and the only ($n-k$)-set which may be missing from $\A$ is $[n]-\Gamma(u)$. However, if $\A$ contains all the $k$-sets then it contains this set as well, so we are done.

So assume $uv \notin E(G)$, and that each of $u$ and $v$ is missing $\ds\frac{n-1}{2}$ edges. Then the only ($n-k$)-sets which may be missing are $[n]-\Gamma(u)$ and $[n]-\Gamma(v)$. As before if either of these is actually in $\A$ then we are done. But if $\A$ contains any $k$-set which is a subset of one of these ($n-k$)-sets and contains both $u$ and $v$, then $\A$ must also contain that ($n-k$)-set, and since there are at least two such $k$-sets when $n \ge 7$, and the remaining cases are trivial, we are done.
\end{proof}

\section{$FC(k)$-families for small values of $k$}\label{smallsec}

First let us consider the case $k = 5$. In \cite{Vaug}, Vaughan showed that a $5$-set with all its $4$-subsets, and a $5$-set with four of its $4$-subsets and four of its $3$-subsets are $FC(5$)-families, using Poonen's Theorem with $\ul{c} = 1$. By using different values of $c_i$ and a (fairly simple-minded) computer program, we have been able to show much more, characterising exactly the $FC(5)$-families.

\begin{proof}[Proof of Theorem~\ref{fvs}]
First we show that the given families are $FC$, using Poonen's Theorem. The required inequalities, $K(\A) \ge 0$ for all $\A \subset \P(5)$ such that $\A \uplus \B \subseteq \A$, follow by a tedious case analysis (either by hand or by computer) once the correct values of $c_i$ have been identified. The search is narrowed by considering only solutions to the five inequalities given by $\A = \B \uplus \P([5]\setminus \{i\})$ for $i = 1, \ldots ,5$.

The following are examples of $\ul{c}$'s which work: \\ (1) If the three 3-sets are contained in some 4-set then let $\ul{c} = 1$ (this is Corollary 4 of \cite{Poon}). If the 3-sets cover [5] then there are three cases to consider:\\[-2ex]

\hspace{1cm} 1. $\{1,2,3\}$, $\{1,2,4\}$, $\{1,2,5\}$: $\ul{c} = (3,3,2,2,2)$ will do;

\hspace{1cm} 2. $\{1,2,3\}$, $\{1,2,4\}$, $\{1,3,5\}$: $\ul{c} = (6,5,5,3,3)$ works;

\hspace{1cm} 3. $\{1,2,3\}$, $\{1,2,4\}$, $\{3,4,5\}$: $\ul{c}= (2,2,2,2,1)$ is an example. \\[+1ex] (It is easy to see that all other possibilities are just permutations of one of these three.) \\
For the rest of the set systems listed in Theorem~\ref{fvs}, suitable values of $\ul{c}$ are\\[+0.5ex]
(2) $\ul{c} = (24,22,19,19,4)$ \\[+0.5ex]
(3) $\ul{c} = (4,3,3,3,3)$\\[+0.5ex]
(4) $\ul{c} = (4,3,3,3,3)$\\[+0.5ex]
(5) $\ul{c} = (14,14,14,9,9)$\\[+0.5ex]
(6) $\ul{c} = 1$ (This is Theorem 4.2 of \cite{Vaug}.)\\

Now we show that these are the only proper $FC(5)$-families (up to permutations of $\{1,2,3,4,5\}$). It suffices to show that the following families are not $FC$.\\[-2ex]

\hspace{1.1cm} 1. $\{1,2,3\}$, $\{1,4,5\}$, $\{1,2,3,4\}$, $\{1,2,3,5\}$ and $\{1,2,4,5\}$

\hspace{1.1cm} 2. $\{1,2,3\}$, $\{1,2,4\}$, $\{1,2,3,4\}$, $\{1,2,3,5\}$ and $\{1,2,4,5\}$

\hspace{1.1cm} 3. $\{1,2,3\}$, $\{1,2,3,4\}$, $\{1,2,3,5\}$, $\{1,2,4,5\}$ and $\{1,3,4,5\}$\\[-2ex]

In each case we show that no values of $c_i$ simultaneously satisfy the inequalities given by $K(\A) \ge 0$, with $\A = \B \uplus \P([5]-j)$ for $j = 1, \ldots ,5$. The result then follows by Poonen's Theorem.

1. Noting that $c_1 + c_2 + c_3 + c_4 + c_5 = 1$, the five inequalities reduce to
\begin{eqnarray}
c_1 & \le & \ds\frac{1}{4} \\
4c_1 - 12c_2 + 2c_3 & \ge & -1\\
9c_1 + 5c_2 + 5c_3 & \ge & 4 \\
c_1 + c_2 & \ge & 3c_3
\end{eqnarray}

From (1) and (2) we have $c_2 \le \ds\frac{c_3 + 1}{6}$ and from (1) and (3) we have $c_2 + c_3 \ge \ds\frac{7}{20}$. But now, since $\ds\frac{7}{20} \le c_2 + c_3 \le \ds\frac{7c_3 + 1}{6} \Rightarrow c_3 \ge \ds\frac{11}{70}$, $3c_3 \le c_1 + c_2 \le \ds\frac{2c_3 + 5}{12} \Rightarrow c_3 \le \ds\frac{5}{34}$ and $\ds\frac{5}{34} < \ds\frac{11}{70}$, we have a contradiction.

2. Form the inequalities as before, and noting the symmetries in $\B$, let $x = c_1 + c_2$, $y = c_3 + c_4$ and $z = c_5 + c_6$. Adding inequalities and making the substitution $z = 1 - x - y$ gives
\begin{eqnarray}
x & \le & y \\
2x & \ge & 3y\\
16x + 14y & \ge & 13
\end{eqnarray}
But (5) and (6) $\Rightarrow x = y = 0$, so again we have a
contradiction.

3. Again noting the symmetries we let $x = c_1$ and $y = c_2 + c_3$, and reduce as before to get
\begin{eqnarray}
-6x + y & \ge & -1 \\
4x - 5y & \ge & -1  \\
9x + 7y & \ge & 5
\end{eqnarray}
But (8) and (9) $\Rightarrow x \le \ds\frac {3}{13}$ and $y \le \ds\frac{5}{13}$, so $9x + 7y \le \ds\frac{62}{13} < 5$, contradicting (10).
\end{proof}

We next consider $FC(6)$-families. Vaughan~\cite{Vaug} proved that any ten 4-sets, or any eight 4-sets together with six 5-sets, generate an $FC$-family. Theorem~\ref{sxs} improves these results, and also gives the exact number of 3-sets in $\{1, \ldots ,6\}$ which are guaranteed to generate an $FC$-family.

\begin{proof}[Proof of the first half of Theorem~\ref{sxs}]
The upper bound again follows by computer-based case analysis. If any three of the $3$-sets do not cover all six elements then we are done by Theorem~\ref{fvs}. Hence each element of $[6]$ must lie in at least (and so exactly) two of the 3-sets. If some pair of 3-sets intersect in two elements then the 3-sets must be $\{1,2,3\}$, $\{1,2,4\}$, $\{3,5,6\}$ and $\{4,5,6\}$ (under some permutation of $\{1,2,3,4,5\}$). Otherwise the sets are $\{1,2,3\}$, $\{1,4,5\}$, $\{2,4,6\}$ and $\{3,5,6\}$. In either case condition 2 of Poonen's Theorem is satisfied by $\ul{c} = 1$.

The lower bound follows by applying the method of Theorem~\ref{fvs} to the family generated by $\{1,2,3\}$, $\{1,2,4\}$, and $\{3,5,6\}$. Let $x = c_1 + c_2$, $y = c_3$ and $z = c_4$, and get
\begin{eqnarray*}
y + z & \ge & x \\
4y + z & \le & 1 \\
x & \ge & 3z \\
x + 2y + z & \ge & 1
\end{eqnarray*}

Now, from the first and third equations we get $y \ge 2z$, and so by the second, $z \le \ds\frac{1}{9}$. Hence $x + 2y + z \le 3y + 2z \le 3\left( \ds\frac{1-z}{4} \right) + 2z = \ds\frac{3}{4} + \ds\frac{5z}{4} < 1$, a contradiction.
\end{proof}

Unfortunately from now on the case analysis involved in proving upper bounds on $FC(k,n)$ becomes too lengthy to be performed by our computer program, so for the second half of Theorem~\ref{sxs} we need to be a little more clever.

\begin{proof}[Proof of the second half of Theorem~\ref{sxs}]
Let $\B$ be generated by eight 4-sets, and note that each element of $[6] = \{1,2,3,4,5,6\}$ must be contained in at least four of the 4-sets, else we are done by Theorem 4.2 of \cite{Vaug}. We shall show that in all but a few special cases the result follows by Poonen's Theorem with $\ul{c} = 1$, and then deal with those cases separately.

Recall that $n_i$ denotes the number of $i$-sets in $\A$. In all cases we may assume that $n_1 \le 5$, and by Theorem~\ref{vcj} we may assume that $n_2 \le 10$ when $\ul{c} = 1$. Also recall Lemma A and observe that a trivial calculation (as in \cite{Vaug}) gives us $d_{2}(6,4,r) \le 10 - r$ for all $1 \le r \le 10$.

Suppose first that the eight 4-sets generate at least 5 of the 5-sets, and let $\ul{c} = 1$, so $K(\A) = \ds\sum_{i=1}^3 i \: (n_{3+i}-n_{3-i})$. Then $n_5 \ge n_1$, since if $\A$ is missing at least two of the 5-sets in $\B$ then by Lemma A $n_1 \le d_{1}(6,5,2) = 0$, if $\A$ is missing exactly one of the 5-sets in $\B$ then $n_1 \le d_{1}(6,5,1) = 1$, and if $\A$ contains all of the 5-sets in $\B$ then $n_5 \ge 5 \ge n_1$ by assumption. Also, since $d_{2}(6,4,r) \le 10 - r$ for all $1 \le r \le 8$, Theorem~\ref{vcj} and Lemma A give us $n_4 - n_2 \ge -2$.

Now, if $\emptyset \notin \A$ then $K(\A) \ge 0$ , so assume that $\emptyset\in \A$. Then $n_5 \ge 5$ and $n_4 \ge 8$, and so if $K(\A) < 0$, we must have $n_1 = n_5 = 5$. But the 1-sets now generate at least one new 4-set, so $n_4 \ge 9$ and $n_2 = 10$. We are done if each 5-set in [6] contains at most three 4-sets of $\B$, since then $n_4 \ge 10$, so assume some four of the 4-sets lie in $\{1,2,3,4,5\}$. Also note that since $n_5 = 5$, at least seven of the 4-sets in $\B$ contain $i$, where $[6]\setminus \{i\}$ is the missing 5-set.\\

\noindent\ul{Claim 1}: If at least seven of the 4-sets contain 1, no other element is contained in more than six, and four of the 4-sets are contained in $\{1,2,3,4,5\}$, then the conditions of Poonen's Theorem hold with $\ul{c} = (11,9,9,9,9,9)$.

\begin{proof}[Proof of claim]
The values of $c_i$ were chosen by considering $\P([6]\setminus \{i\}) \uplus \B $ for $i = 2, \ldots ,6$, assuming $c_2 = \ldots = c_6$ and minimising $c_1$. Let $\A \subset \P(5)$ such that $\A \uplus \B \subseteq \A$, and note that we still have that $\A$ contains at least $n_1$ of the 5-sets and $n_2 - 2$ of the 4-sets of $\B$. Recall that $N_i$ denotes the contribution to $K(\A)$ of the $i$-sets in $\A$, and note that $N_3 \ge -20$.

Consider first the case $\emptyset \notin \A$. Then $N_6 + N_0 = +56$, since we may assume $[6]\in \A$. Also $N_4 + N_2 \ge -44$, since 2-sets contribute at least $-20$, and all but at most one of the 4-sets in $\B$ contributes 20, the other contributing 16. Suppose $\{2,3,4,5,6\} \notin \A$. Then we have $N_3 \ge -8$, and $N_5 + N_1 \ge 0$, so $K(\A) \ge 0$. If $\{2,3,4,5,6\}\in \A$ then $N_5 + N_1 \ge 34$ (assuming $n_1 \le 5$), so since $N_3 \ge -20$, $K(\A) \ge 0$ still holds.

So assume now that $\emptyset\in \A$, so $n_4 \ge 8$ and $n_5 \ge 5$ by Observation~\ref{triv}. First note that if $n_1 \le 3$ then $N_5 + N_1 \ge 76$, $N_4 + N_3 + N_2 \ge -64$ and we are done. So assume $n_1 \ge 4$ and again suppose that $\{2,3,4,5,6\} \notin \A$. If $n_1 = 4$ then $N_5 + N_1 \ge 38$, $N_3 \ge -8$ and so if $K(\A) < 0$ then $n_2 = 10$ and $n_4 = 8$. But this is impossible since any four 1-sets and ten 2-sets clearly generate either $\{2,3,4,5,6\}$ or a 4-set not in $\B$. If $n_1 = 5$, then since $\{2,3,4,5,6\} \notin \A$, the only possibilities for $\A$ are $\P([6]\setminus \{i\}) \uplus \B$ and $\P([6]\setminus \{i\}) \uplus \B \uplus \{1,i\}$, where $i \in \{2,3,4,5,6\}$ (we may assume no extra 3-sets containing 1 are added to $\A$, since these increase $K(\A)$). We chose $\ul{c}$ so that $K(\A) \ge 0$ in the former case, and the latter gives a higher value of $K(\A)$, since new 4-sets are generated.

Hence we may assume that $\{2,3,4,5,6\} \in \A$, and so $n_5 = 6$. If $n_1 \le 4$ we get $N_5 + N_1 \ge 72$ and are done as before, so assume also that $n_1 = 5$. If $\{1\} \notin \A$ we get $N_5 + N_1 = 34$ and $n_4 \ge 12$, but $n_2 \ge 13 \Rightarrow n_4 = 15$, so we are done. But if $\{1\}\in \A$ we get $N_5 + N_1 = 38$, $N_3 \ge -8$ and $N_4 + N_2 \ge -24$ (since $n_4 \ge 9$), so again $K(\A) \ge 0$, and the proof of claim 1 is complete.
\end{proof}

So from now on we may assume that $\B$ contains at most four 5-sets. But since any two 4-sets in a 5-set generate that 5-set (so a missing 5-set implies four missing 4-sets), and each 4-set is contained in only two 5-sets of [6], eight 4-sets in [6] must generate at least four 5-sets. Hence $\B$ contains exactly four 5-sets. But now we can reduce the problem to a single family, for suppose wlog that $\{1,3,4,5,6\}$ and $\{2,3,4,5,6\}$ are the 5-sets missing from $\B$. Then $\B$ must contain all 4-sets which contain both 1 and 2, one containing 1 and not 2, and one containing 2 and not 1. These last two sets must have intersection 2, as otherwise we would have five 4-sets in $\{1,2,3,4,5\}$, and by symmetry we may take any pair with this property.

It follows that the proof of the upper bound is completed by the following claim.\\

\noindent\ul{Claim 2}: If $\B$ is generated by the sets $\{1,2,3,4\}$, $\{1,2,3,5\}$, $\{1,2,3,6\}$, $\{1,2,4,5\}$, $\{1,2,4,6\}$, $\{1,2,5,6\}$, $\{1,3,4,5\}$ and $\{2,3,4,6\}$, $\A \subset \P(5)$, $\A \uplus \B \subset \A$ and $\ul{c} = (8,8,7,7,7,7)$, then $K(\A) \ge 0$.

\begin{proof}[Proof of claim]
As in Claim 1, the values of $c_i$ were chosen by looking at $\P([6]\setminus \{i\}) \uplus \B$ for $1 \le i \le 6$, letting $c_3 = \ldots = c_6$ and minimising $c_1 + c_2$. Our approach follows the same lines as the proof of claim 1. Note first that if $n_1 = 5$, either $\{1,3,4,5,6\} \in \A$ or $\{2,3,4,5,6\} \in \A$, so $n_5 \ge 5$ and as before $\A$ still contains at least $n_2 - 2$ of the 4-sets and $n_1$ of the 5-sets of $\B$.

Consider first the case $\emptyset \notin \A$. We have $N_6 + N_0 = 44$, $N_5 + N_1 \ge 0$, $N_4 + N_2 \ge -36$ (since 2-sets contribute at least $-16$, and six of the 4-sets in $\B$ contribute 16, the other two contributing 14), and $N_3 \ge -4$, so $K(\A) \ge 0$.

So we may assume that $\emptyset\in \A$, so $n_4 \ge 8$ and $n_5 \ge 4$ by Observation~\ref{triv}. If $n_5 = 4$, then $\{5\}, \{6\} \notin \A$ and $n_2 \le 8$, since the only 2-sets containing 5 or 6 which can be in $\A$ are $\{1,5\}$ and $\{2,6\}$. Hence $N_5 + N_1 \ge 4$ and $N_4 + N_2 \ge (16 \times 6) + (14 \times 2) - 16 - (14 \times 6) - 12 = 12$. We still have $N_3 \ge -4$, so $K(\A) \ge 0$.

Suppose next that $n_5 = 5$, and $\{1,3,4,5,6\} \in \A$, say. If $n_1 \le 3$ then $N_5 + N_3 + N_1 \ge 54$ and we are done. If $n_1 = 4$ then $N_5 + N_3 + N_1 \ge 26$, so $n_2 = 10$ and $n_4 = 8$, which is impossible. Hence $n_1 = 5$, and $\{5\} \notin \A$, so the only possibilities for $\A$ are $\P([6]-5) \uplus \B$ and $\P([6]-5) \uplus \B \uplus \{1,5\}$ (since $\{2,3,4,5,6\} \notin \A$, and we may assume no extra sets with non-negative contribution to $K(\A)$ are added). But these both give $K(\A) \ge 0$, so we are done with this case also.

So we may assume that $n_5 = 6$. Now if $n_1 \le 4$ then $N_5 + N_3 + N_1 \ge 52$, and so we are done. But if $n_1 = 5$ then $N_5 + N_3 + N_1 \ge 24$ and $n_4 \ge 10$, since no 5-set in [6] contains more than three 4-sets of $\B$, so $N_4 + N_2 \ge -12$, and the claim follows.
\end{proof}

For the lower bound consider the set system $\S = \{\{1,2,i,j\} : \{i,j\} \subset \{3,4,5,6\}\}$ and let $\B$ be generated by $\S$. We apply the usual method. Let $x = c_1 + c_2$ and $y = c_3 + c_4 + c_5 + c_6$, and add the inequalities to get
\begin{eqnarray*}
38x + 46y & \ge & 43 \\
92x + 67y & \ge & 78
\end{eqnarray*}
which imply that $x \le \ds\frac{1}{4}$ and $x \ge \ds\frac{11}{25}$ respectively. But $\ds\frac{1}{4} < \ds\frac{11}{25}$, so we have a contradiction. This implies that $FC(4,6) \ge 7$.
\end{proof}

We shall also use the following two simple results in Section 5.

\begin{lemma}\label{svnthr}
$FC(3,7) \le 6$.
\end{lemma}

\begin{proof}
Consider any six 3-sets in $\{1,2,3,4,5,6,7\}$. Some element is contained in at most two of them, so we have at least four 3-sets composed solely of the other six elements. The result now follows by Theorem~\ref{sxs}.
\end{proof}

\begin{lemma}\label{svnfr}
$FC(4,7) \le 18$.
\end{lemma}

\begin{proof}
Suppose we have eighteen 4-sets in $\{1,2,3,4,5,6,7\}$. Then some
element must be contained in no more than ten of them, so we have
at least eight 4-sets contained in a 6-set, and the result
follows by Theorem~\ref{sxs}.
\end{proof}

\section{The case $n = 9$ of the conjecture}\label{9sec}

We now provide an application of the above results by proving the Union-Closed Sets Conjecture in the case that the size of the largest set is at most nine. This improves the previous known bound by one. The idea of the proof is to show that if the family contains none of the above $FC$-families, then either the average size of a set is at least $\ds\frac{1}{2}$, or the family contains very few sets. We shall need the following results for the proof. The first two were proved by Gao and Yu~\cite{GY}, the third by Poonen~\cite{Poon}.

\begin{thmC}\label{mlim}
If $m \le 32$ then the UC-Sets Conjecture holds for $\A$.
\end{thmC}

\begin{thmD}\label{eight}
If $n \le 8$ then the UC-Sets Conjecture holds for $\A$.
\end{thmD}

\begin{lemmaE}\label{big}
We may assume that $\A$ contains at least two $(n-1)$-sets.
\end{lemmaE}

\begin{lemma}\label{lots}
One $3$-set and thirteen $4$-sets in $[9]$ generate either at least two $7$-sets, at least three $6$-sets, or an $FC$-family.
\end{lemma}

\begin{proof}
Suppose only at most one 7-set and two 6-sets are generated. Partition the 4-sets according to the size of their intersection with the 3-set. If we have more than six 4-sets with intersection 2, then we get at least four distinct 2-sets by removing the elements of our 3-set from them (by Theorem~\ref{fvs} (5)), which generate either at least two 4-sets or three 3-sets, a contradiction. But we can have at most three 4-sets with intersection 3, at most two with intersection 1 and at most one with intersection 0.
\end{proof}

\begin{lemma}\label{asvn}
Five $3$-sets in $[9]$ generate either an $FC$-family or a $7$-set.
\end{lemma}

\begin{proof}
Assume not and choose three of the 3-sets in such a way as to maximise the size of their union. Let this maximum be $t$. By Theorem~\ref{fvs}, $t \ge 6$.\\[+2ex]
Case 1: $t = 9$. The addition of any other 3-set forms a 7-set so we are done.\\[+2ex]
Case 2: $t = 8$. Suppose without loss that the 3-sets are $\{1,2,3\}$, $\{4,5,6\}$ and $\{6,7,8\}$. Then by Theorem~\ref{sxs} we must add a 3-set which intersects $\{1,2,3\}$, and the only (type of) 3-set which does so and whose addition does not create a 7-set is $\{1,4,5\}$. Now, consider the possible intersections of the fifth 3-set with $\{1,2,3,4,5\}$. By Theorem~\ref{fvs} it cannot have order 3; by Theorem~\ref{sxs} if it has order 2 then we have a 7-set; if it has order 1 we clearly have a 7-set, so the fifth 3-set must be contained in $\{6,7,8,9\}$. But then it must contain 9, so we have the 7-set $\{1,4,5,6,7,8,9\}$.\\[+2ex]
Case 3: $t = 6$. Each of the other 3-sets can contain at most one of these six elements (if two we'd be able to form a 7-set, if three we'd have an $FC$-family by Theorem~\ref{sxs}), so by maximality no two of the chosen 3-sets intersect in two elements. Hence wlog the chosen 3-sets are $\{1,2,3\}$, $\{1,4,5\}$ and $\{2,4,6\}$. At least one of the remaining sets must be of the form $\{i,7,8\}$ say, with $i \in \{1,2,3,4,5,6\}$, but then we have a 7-set.
\end{proof}

\begin{proof}[Proof of Theorem~\ref{nine}]
Suppose $\A$ is a union-closed family in $\P(9)$ for which the conjecture fails, and observe that by Theorem~D we may assume that $[9] \in \A$. Also note that we may assume that $n_1 = n_2 = 0$, since every 1- and 2-set forms an $FC$-family. We are trying to show $K(\A) \ge 0$ for $\ul{c} = 1$, so an $r$-set contributes $2r - 9$ to $K(\A)$. We shall consider the contribution of the (non-empty) sets in the 7-sets of [9]. Notice that there are $\ds{9 \choose 7} = 36$ of them, and that the contribution of [9], $\emptyset$ and the 8-sets of $\A$ is at least $+14$.

Consider a single 7-set of [9], not necessarily in $\A$. If the 7-set itself is in $\A$ then it contributes +5 to $K(\A)$, but each $r$-subset lies in $\ds{{9-r} \choose {7-r}}$ different 7-sets, so we must divide its contribution by this number. Summing over all 7-sets will then give us the total contribution of the 3-, 4-, 5-, 6- and 7-sets to $K(\A)$. Hence 6-sets contribute +1, 5-sets $+\ds\frac{1}{6}$, 4-sets $-\ds\frac{1}{10}$ and 3-sets $-\ds\frac{1}{5}$.

We need a lower bound on the contribution of a given 7-set's non-empty subsets. Let $r_i$ denote the number of $i$-sets in this 7-set, and divide into cases as follows:\\[+2ex]
Case 1: The 7-set is in $\A$. Then the contribution is at least $5 - \ds\frac{r_4}{10} - \ds\frac{r_3}{5} \ge 5 - \ds\frac{17}{10} - \ds\frac{5}{5} = \ds\frac{23}{10}$, by Lemmas~\ref{svnthr} and \ref{svnfr}.\\[+2ex]
Case 2: The largest subset in $\A$ has six elements. Then the contribution of the subsets is at least 0, since $r_4 \le 7$ and $r_3 \le 3$ by Theorem~\ref{sxs}, so we may assume that either $r_3 = 2$ and $r_4 = 7$, or $r_3 = 3$. But then $r_5 \ge 1$, so we may assume that $r_3 = 3$ and $r_4 \ge 6$. But now $r_5 \ge 2$ and we are done.\\[+2ex]
Case 3: The largest subset in $\A$ has five elements. Then it follows from Theorem~\ref{fvs} that the worst case is $r_3 = 2$, $r_4 = 3$, which gives a total contribution of $-\ds\frac{16}{30}$.\\[+2ex]
Case 4: The largest subset in $\A$ has fewer than five elements. Then $r_3 \le 2$ and $r_4 \le 1$, so the contribution is at least $-\ds\frac{1}{2}$.\\[-2ex]

\noindent\ul{Claim}: If $\A$ contains a 7-set then $K(\A) \ge 0$.

\begin{proof}[Proof of claim]
First note that if $n_8 \ge 3$ then we are done, since then $K(\A) \ge 21 - (36 \times -\ds\frac{16}{30}) > 0$, so assume that $n_8 = 2$. Suppose that $n_7 \ge 2$. Then $K(\A) \ge 14 + 4.6 - (34 \times -\ds\frac{16}{30}) > 0$, so assume that $n_7 = 1$ and observe that if $K(\A) < 0$ then at least 31 of the 7-sets contain no 6- or 7-set, since $14 + 2.3 - (30 \times -\ds\frac{16}{30}) > 0$. Hence $n_3 \le \ds\frac{(31 \times 2) + (4 \times 3) + 5}{15} < 6$ since each 3-set appears in 15 different 7-sets, and $n_4 \le \ds\frac{(31 \times 4) + (4 \times 7) + 17}{10} < 17$, since each 4-set appears in 10 different 7-sets. It follows that the total contribution to $K(\A)$ of all except the 5- and 6-sets is at least $14 + 5 - 16 - 15 = -12$, so assume that the 5- and 6-sets contribute at most $+11$. Note also that $\ds\sum_{i \neq 5,6} n_i \le 26$.

Now, from above if $n_3 = 0$ then $K(\A) > 0$, so by Lemma~\ref{lots}, either $n_6 \ge 3$, or $n_4 \le 12$. In the former case we get $n_5 + n_6 \le 5$, so $m \le 31$ and we are done by Theorem~C. In the latter case the 5- and 6-sets can contribute at most $+7$ and $\ds\sum_{i \neq 5,6} n_i \le 22$, which gives $m \le 29$ and we are again done.
\end{proof}

So we may assume that $n_7 = 0$, and by Lemma~\ref{asvn} we have $n_3 \le 4$. Now, let $p = |\{$7-sets with $r_4 \le 3\}|$ and $q = |\{$7-sets with $r_4 = 4\}|$. Then Theorem~\ref{fvs} and the observations above imply that $\ds\frac{16}{30}p + \ds\frac{13}{30}q > 14$ if $K(\A) < 0$. But then $10n_4 \le 3p + 4q + 7(36 - p - q) = 252 - (4p + 3q) \le 252 - 105 + \ds\frac{q}{4} \le 156$, so $n_4 \le 15$.

Hence the total contribution to $K(\A)$ of all except the 5- and 6-sets is at least $14 - 15 - 12 = -13$, so we may assume that the 5- and 6-sets contribute at most $+12$, and we have $\ds\sum_{i \neq 5,6} n_i \le 23$. The result now follows exactly as before, using Lemma~\ref{lots} and Theorem~C: $n_3 \ge 1$ so either $n_6 \ge 3$, in which case $n_5 + n_6 \le 6$ and $m \le 29$; or $n_4 \le 12$, in which case $n_5 + n_6 \le 9$, and $m \le 29$.
\end{proof}

Although our method does not seem to easily extend to larger values of $n$ (at least, not without first improving our upper bounds on $FC(k,n)$), it does give a short proof of Theorem~D.

\begin{proof}[Alternative proof of Theorem~D]
Suppose $\A$ is a counter-example and let $b = max_{A \in \A} |A|$. Since $n_1 = n_2 = 0$, we may assume that $b \ge 7$.

Case 1: $b = 7$. By Lemma~\ref{svnthr} we have $n_3 \le 5$, and by Lemma E $n_6 \ge 2$, so the average set size is at least $\ds\frac{34}{9} > \ds\frac{7}{2}$ and we are done.

Case 2: $b = 8$. Apply the method of the proof of Theorem ~\ref{nine}. There are $\ds{8 \choose 6} = 28$ 6-sets, each contributing at least $-\ds\frac{1}{5}$ to $K(\A)$ (two 3-sets in a 4-set). But now $K(\A) \ge (2 \times 3) - \ds\frac{28}{5} > 0$, so again we are done.
\end{proof}

We also get the following improvement without doing any extra work.

\begin{thm}\label{newm}
The Union-Closed Sets Conjecture holds in the case $m \le 36$.
\end{thm}

\begin{proof}
Follow the exact method of \cite{GY} and \cite{Poon}, replacing the bound $n \ge 9$ with $n \ge 10$.
\end{proof}

\begin{prob}
Can we improve Theorem~$\ref{mlim}$ further? A useful step would be (using the notation of \textup{\cite{GY}}) to get good lower bounds on $S_r$ in terms of $S_1$, given a condition limiting the average size of the sets in $\A$.
\end{prob}

\begin{prob}
Can we generalise or improve Lemmas$~\ref{lots}$ and $\ref{asvn}$? In particular, how many $k$-sets in $[n]$ do we need to guarantee either an $FC$-family or an $r$-set?
\end{prob}

\section{General bounds on $FC(k,n)$}\label{bigsec}

We have found some values of $FC(k,n)$ for small $k$ and $n$, but if we hope to use our method to solve the conjecture we shall need good asymptotic bounds. Good upper bounds seem hard to prove, but by Proposition 1.4 of \cite{GY} the function is at least defined for any $k$ for sufficiently large $n$. We give the following short proof of the result.

\begin{thm}
For any $k \ge 1$ and $n \ge 2k - 2$, the family $\B$ generated by all the $k$-sets in $[n]$ is an $FC$-family, and hence $FC(k,n) \le \ds{n \choose k}$.
\end{thm}

\begin{proof}
The result follows by Poonen's Theorem, and the following claim.\\[-2ex]

\noindent\ul{Claim}: If $\A \subset \P(n)$ and $\A \uplus \B \subseteq \A$, then $n_{n-r} \ge n_r$ for all $r \le \left\lfloor \ds\frac{n}{2} \right\rfloor$.
\begin{proof}[Proof of claim]
Note that $\B \supset \{B \in \P(n) : |B| \ge k\}$. So if an $r$-set $A$ is in $\A$, then all those $(n-r)$-sets in $\P(n)$ that contain $A$ are also in $\A$. Hence each $r$-set in $\A$ generates $\ds{{n-r} \choose {n-2r}}$ $(n-r)$-sets. Conversely each $(n-r)$-set in $\A$ is generated by at most $\ds{{n-r} \choose r} = {{n-r} \choose {n-2r}}$ $r$-sets, so we are done.
\end{proof}
Now let $\ul{c} = 1$, so $K(\A) \ge 0$ for families satisfying the conditions of the claim. The result follows by Poonen's Theorem.
\end{proof}

\begin{rmk}
This result is almost certainly not best possible, but as yet we have been unable to improve it. The usual lower bound method gives that $\B$ is not $FC$ if $\: (n-2)\ds{{n-2} \choose {k-2}} < 2\ds\sum_{i=0}^{k-3} {{n-2} \choose i}$, which holds if and only if \: $k \ge \ds\frac{n}{2} + \left( \ds\frac{1}{2\sqrt{2}} + o(1) \right) \sqrt{n \operatorname{log}n}$.
\end{rmk}

We have the following bounds on $FC(k,n)$ (the first upper bound, $FC(3,n) \le \ds\frac{2n}{3}$, is due to Vaughan~\cite{Vthrees}). The proofs of the lower bounds are not difficult but lengthy, so we omit them, giving only the main ideas.

\begin{thm}
\hspace{0cm}
\begin{enumerate}
\item[$(a)$] $\left\lfloor \ds\frac{n}{2} \right\rfloor + 1 \le FC(3,n) \le
\ds\frac{2n}{3}$,
\item[$(b)$] $\left(\ds\frac{6}{25} + o(1)\right)n^2 \le FC(4,n) \le
\ds\frac{7}{360}n^4$,
\item[$(c)$] $c_kn^{k-2} \le FC(k,n) \le ex(n,K_{2k-2}^{(k)}) \le \ds{n \choose
k}$ with $c_k > 0$ constant.
\end{enumerate}
\end{thm}

\begin{proof}[Sketch of proof]
The upper bounds are easily obtained by showing that any family in $\P(n)$ with the given number of $k$-sets contains one of the families we (or Vaughan) have already shown to be $FC$. For the lower bounds, first let
\begin{align*}
\B(n,k,r) \; = & \; \{ B \in \P(n) : |B| = k\textup{ and for some }0 \le i \le r - 1\textup{, either }\{4i+1,\\
& \; \; 4i+2, 4i+3\}\textup{ or }\{ 4i+1, 4i+2, 4i+4\}\textup{ is the initial segment of }B \}.
\end{align*}
The bounds are obtained by applying the method of the proofs of Theorems~\ref{fvs} and \ref{sxs}, using the families $\P([n]\setminus \{i\})$ for $1 \le i \le n$, to

(1) $\B = \B(n,3,r)$, if $n = 4r$ or $4r + 1$, or

\hspace{0.5cm} $\B = \B(n,3,r) \cup \{4r-1,4r+1,4r+2\}$ if $n = 4r + 2$ or $4r + 3$.

(2) $\B = \B(n,4,r)$, where $r = \left\lfloor \left( \ds\frac{1}{5} - \eps \right)n \right\rfloor$ for some $\eps > 0$.

(3) $\B = \B(n, k, c_k^{\prime} n)$ for some sufficiently small $c_k^{\prime} > 0$.
\end{proof}

\section{Questions and Conjectures}\label{QCsec}

Several avenues for further research spring readily to mind. For example, given our experience so far we might hope to show that in Poonen's Theorem it is sufficient to consider only the families $\P([n]\setminus \{i\})$ for $1 \le i \le n$. Unfortunately it is not true that these families give exactly the permissible values of $\ul{c}$, as the following example shows.

\begin{counter}
Let $\B = \{ \{1,2,3\}, \{1,2,4\}, \{3,4,5\} \}$, and let $\ul{c} = (9,7,12,12,8)$. Then $K(\A) \ge 0$ for $\A = \P([5]\setminus \{i\}) \uplus \B$ for each $i$, but $K(\A) < 0$ for $\P(\{2,3,4,5\}) \uplus \{1,2\} \uplus \B$.
\end{counter}

However, it is still possible that the following question has an affirmative answer.

\begin{qu}
Do the inequalities $K(\A) \ge 0$ given by the families $\P([n] \setminus \{i\})$ for $1 \le i \le n$ permit \textit{some} solution $\ul{c}$ only when one is possible for all $\A \subset \P(n)$ such that $\A \uplus \B \subseteq \A$?
\end{qu}

In any case we are still inclined to believe the following conjectures, the first of which was suggested (though not specifically conjectured) in \cite{Vthrees}.

\begin{conj}
$FC(3,n) = \left\lfloor \ds\frac{n}{2} \right\rfloor + 1$ for all $n \ge 4$.
\end{conj}

\begin{conj}
$FC(k,n) = \Theta(n^{k-2})$ for all $k \ge 2$.
\end{conj}

A final conjecture is suggested simply on the basis that it seems plausible, and is true for the $FC$-families we have discovered so far.

\begin{conj}
Suppose a union-closed family $\B$ has minimal generating family $\S$, and let $\B^\prime$ be generated by $(\S \setminus \{B\}) \cup \{B \cup \{i\} \}$ for some set $B \in \S$. If $\B$ is not $FC$, then $\B^\prime$ is not $FC$ either.
\end{conj}

\end{document}